\documentclass[preprint]{elsarticle}
\usepackage{amsmath}
\usepackage{amssymb}
\usepackage{epsfig}
\usepackage{epstopdf}
\usepackage[latin5]{inputenc}
\usepackage{subfigure}
\usepackage{lineno,hyperref}

\modulolinenumbers[5]

\newtheorem{proposition}{Proposition}
\newtheorem{theorem}{Theorem}
\newtheorem{example}{Example}
\newtheorem{remark}{Remark}
\newtheorem{lemma}{Lemma}
\newtheorem{corollary}{Corollary}

\newproof{proof}{Proof}

\usepackage{color}

\numberwithin{equation}{section}

\begin{document}

\title{$K^{\alpha}-$translators of offset surfaces}

\author{Burcu Bekta\c s Demirci}
\ead{bbektas@fsm.edu.tr}

\address{Fatih Sultan Mehmet Vak{\i}f University, Faculty of Engineering, Department of Software Engineering, 34025, Istanbul, T\"{u}rkiye.}

\author{Ferda\u{g} Kahraman Aksoyak}
\ead{ferdag.aksoyak@ahievran.edu.tr}

\address{K{\i}r\c{s}ehir Ahi Evran University, Faculty of Education, Division of Elementary Mathematics Education, 40100, K{\i}r\c{s}ehir, T\"{u}rkiye.}

\author{Murat Babaarslan*}
\ead{murat.babaarslan@bozok.edu.tr}

\cortext[cor1]{Corresponding author}
\address{Yozgat Bozok University, Faculty of Science and Letters, Department of Mathematics, 66100, Yozgat, T\"{u}rkiye.}

\begin{abstract}
In this paper, we study $K^{\alpha}$--translators on parallel surfaces and canal surfaces in 3-dimensional Euclidean space $\mathbb{E}^3$. First, we investigate the condition under which two parallel surfaces can become $K^{\alpha}$--translators moving with the same speed $w$. Then, we examine $K^{\alpha}$--translators on canal surfaces and we show that if a canal surface is $K^{\alpha}$--translator, then it must be a surface of revolution in $\mathbb{E}^3$. We also provide examples for moving a surface of revolution under $K$--flow (Gauss curvature flow) and $K^{-1/2}$--flow (inverse Gauss curvature flow) along a direction $w=(0,0,1)$ and we illustrate such surfaces using Wolfram Mathematica 10.4. Finally, we prove that no $K^{\alpha}$--translators exist on the parallel surface of a rotational surface obtained from a canal surface with the same speed $w$, while the such rotational surfaces itself is a $K^{\alpha}$--translator with speed $w$.
\end{abstract}

\begin{keyword}
Parallel surface, Canal surface, Linear Weingarten surface, $K^{\alpha}-$translator.
\MSC[2010] 53A05
\end{keyword}

\date{\today}

\maketitle

\section{Introduction}
\label{intro}
In 1974, the Gauss curvature flow was first introduced by Firey \cite{Firey} to describe the shape of worn stones.
Tso \cite{Tso} laid the foundation for subsequent studies on the Gauss curvature flows.
Also, Chow \cite{Chow,Chow2} supported a general study of $K^{\alpha}$--flows (see \cite{Fei}).
In the past forty years, the Gauss curvature flows have been developed fully by the different authors.
This interest reason is that it has important results and applications in both physics and mathematics.
Further studies can be found in \cite{Andrews,Aydin,Aydin2,Chen,Choi,Fei,Urbas,Urbas2,Urbas3,Wu} and the references therein.
Then, the definition of Gauss curvature flow is given as follows.

Let $M$ be a strictly convex surface parametrized by a smooth immersion $X:M\rightarrow\mathbb{E}^3$.
The $K^{\alpha}$--flow is a one-parameter family of smooth immersions
$X_t=X(.,t):M\rightarrow\mathbb{R}^3,\;t\in [0,T)$ such that $X_0=X$ and satisfying
\begin{equation}
\label{defKaflow}
\frac{\partial}{\partial t}X(p,t)=-K(p,t)^{\alpha}U(p,t),\;\;(p,t)\in\Sigma\times [0,T),
\end{equation}
where $\alpha$ is a nonzero constant, $U$ and $K$ denote the unit normal vector field and the Gauss curvature, respectively.
The power $\alpha$ of Gauss curvature $K$ can be positive or negative
and the flow of Gauss curvature by the negative power is called inverse Gauss curvature flow, \cite{Urbas3}.
If $\alpha=1$, then the surface $M$ is spherical, \cite{Andrews2, Aydin}.

If the surface moves under the $K^{\alpha}$--flow along a spatial direction $w\in\mathbb{R}^3$,
then $X_0$ satisfies $K^{\alpha}=\lambda\langle U, w \rangle$ for some constant $\lambda$.
The vector $w$ is unit and it is called as the speed of the flow.
We note that, after an expansion, we can choose $\lambda$ as $1$.
A surface $M$ is called a translator by the $K^{\alpha}$--flow with the speed $w$,
or simply a $K^\alpha$--translator, if the following equation holds
\begin{equation}
\label{deftranslator}
K^\alpha=\langle U, w\rangle.
\end{equation}
(see \cite{Aydin}).

Aydin and L{\'o}pez \cite{Aydin} studied $K^{\alpha}$--translators that the surface is defined kinematically as the movement of a curve
by a uniparametric family of rigid motions under the geometric condition in $\mathbb{E}^3$.
Also, they investigated $K^{\alpha}$--translators of rotational surfaces, where
the rotational axis of rotational surfaces is $z$-axis
and the generating curve of the rotational surface is taken as a plane curve with $z=f(x)$.
In \cite{Aydin2}, Aydin and L{\'o}pez found all $K^{\alpha}$--translators of rotational surfaces
with time-like axis, space-like axis and light-like axis in Minkowski 3-space $\mathbb{E}^3_1$.
$K^{\alpha}$--translating solitons under certain hypotheses on the curvatures in $\mathbb{E}^3$ is characterized by Fei et al. \cite{Fei}.

A paralel surface, also known as offset surface, is a surface whose points are at a constant distance along the normal from another surface. These surfaces were studied in many papers and books, see \cite{Craig} and \cite{Eisenhart}. A generalized offset surface is obtained by replacing the constant distance with variable distance function \cite{Georgiev}. These surfaces were first introduced by Brechner in \cite{Brechner}. Offset surfaces are sometimes defined as the envolope of circle or spheres centered on a generator curve or surface (see \cite{Bruce}). Also, canal surfaces are defined as envelope of a nonparameter set of spheres, centered at a spine curve $c(s)$ with radius $r(s)$. If $r(s)$ is a constant, then the canal surfaces are called as tubes. Thus, canal surfaces can be considered as offset surfaces. More recently, parallel and canal surfaces have been studied in Computer Aided Geometric Design (CAGD) literature under the name offset surfaces (see \cite{Barnhill}). The fundamental geometric and algebraic properties for the canal surfaces were discussed in \cite{Xu}.

The study of surfaces on which there is a relationship between the principal curvatures (Weingarten surfaces) were initiated by Julius Weingarten in 1861 and taken up by a lot of authors for different surfaces (see \cite{Dillen}). Weingarten surfaces are also studied in CAGD and shape investigation because of their attractive curvature properties (see \cite{Grant}).

Weingarten and linear Weingarten tubes in $\mathbb{E}^3$ were investigated by \cite{Ro}. As a generalization, Weingarten and linear Weingarten canal surfaces were studied by \cite{Tuncer}. Parallel linear Weingarten surfaces in $\mathbb{E}^3$ were investigated in \cite{Yayli}.

In this paper, we examine how Gauss curvature flows behave in the case of parallel surfaces,
where all points of the surface are equidistant from a given surface and in the case of canal surfaces.
First, we analyse $K^{\alpha}-$translators on parallel surfaces.
Then, we get some conditions for parallel surfaces to become $K^{\alpha}$--translators with the same speed.
Secondly, we study $K^{\alpha}$--translators on canal surfaces and
we find the relation between the tangent vector of center curve of canal surface and the speed of flow.
Also, we prove that if a canal surface is a $K^{\alpha}$--translator, then it becomes a surface of revolution.
We give some examples by using Wolfram Mathematica 10.4.
Finally, we show that canal surfaces and their parallel surfaces cannot evolve under the same $K^{\alpha}$--flow along the same direction $w$.

\section{Preliminaries}
\label{Sec:1}
In this section, we will state some necessary notations, definitions and theorems about the surfaces (especially, the parallel surfaces) in Euclidean 3-space $\mathbb{E}^3$, see \cite{Oprea}.

Let $M$ be a regular surface in $\mathbb{E}^3$ parametrized by $X(s,\theta)$. The regular surface $M$ is always orientable and its unit normal vector is given by
\begin{equation}
   U=\frac{X_s\times X_\theta}{\left \Vert X_s\times X_\theta \right \Vert}.
\end{equation}
We know that $\nabla_{v}U$ tell us how $M$ curves in the $v-$direction, where $\nabla$ is covariant derivative and $v\in T_{p}M$. Therefore, the shape operator (or Weingarten map) of $M$ at $p$ is defined
\begin{equation}
   S_p(v)=-\nabla_{v}U.
\end{equation}

Moreover, the Gauss curvature and mean curvature of $M$ at $p$ are defined by
\begin{equation}
K=\det(S_p)=k_{1}k_{2},
\end{equation}
and
\begin{equation}
H=\frac{1}{2}\textup{trace}(S_p)=\frac{k_{1}+k_{2}}{2},
\end{equation}
where $k_{1}$, $k_{2}$ are the principal curvatures of $M$.

We note that $M$ is said to be flat if $K(p)=0$ for every $p\in M$ and it said to be minimal if $H(p)=0$ for every $p\in M$.

For $\lambda\in\mathbb{R}-\{0\}$, the parallel surface to $M$ at a distance $\lambda$ is denoted by  $M^{\lambda}$ and it is parametrized as
\begin{equation}
  \bar{X}(s,\theta)=X(s,\theta)+\lambda U(s,\theta).
\end{equation}

For the principal directions $X_s$, $X_\theta$ and  principal curvatures $k_{1}$, $k_{2}$, the derivatives of $M^{\lambda}$ are computed by $\bar{X}_{s}=(1-\lambda k_{1})X_{s}$ and $\bar{X}_{\theta}=(1-\lambda k_{2})X_{\theta}$.
By a simple calculation, the unit normal of $M^{\lambda}$ is given by $\bar{U}=\varepsilon U$, where $\varepsilon$ is sign of $(1-\lambda k_{1})(1-\lambda k_{2})$.
By using the Gauss curvature and mean curvature of $M$, we get $(1-\lambda k_{1})(1-\lambda k_{2})=1-2\lambda H+ \lambda^{2}K$.
Also, the Gauss curvature and mean curvature of $M^{\lambda}$ is given by
\begin{equation}
\label{Gaussparallel}
\bar{K}=\frac{K}{1-2\lambda H+ \lambda^{2}K}
\end{equation}
and
\begin{equation}
\label{meancurvparallel}
\bar{H}=\frac{H-\lambda K}{1-2\lambda H+ \lambda^{2}K}.
\end{equation}

A surface $M$ is called a Weingarten surface if there is a relation between its two principal curvatures $k_{1}$ and $k_{2}$, that is,
if there is a smooth function $V$ of two variables such that $V(k_{1}, k_{2}) = 0$ implies a relation $W(K, H) = 0$, \cite{Yayli}.
If there is a linear relation between the Gauss curvature $K$ and mean curvature $H$ of $M$, that is, $aK+bH=c$, where $a,b,c \in \mathbb{R}, (a,b,c)\neq (0,0,0)$, then $M$ is called a  $(K,H)$--linear Weingarten surface.

We note that when $b=0$, a Weingarten surface $M$ reduces to a surface with constant Gaussian curvature.
Also, when $a=0$, a Weingarten surface $M$ reduces to a surface with constant mean curvature.
Thus, we can say that linear Weingarten surfaces is a natural generalization of surfaces with
constant Gaussian curvature or surfaces with constant mean curvature, \cite{Ro}.

The parallel surfaces of linear Weingarten surfaces in $\mathbb{E}^3$ are linear Weingarten surfaces. Thus, we can give the following lemma.

\begin{lemma}
\cite{Jeromin}
Let  $M$ be a linear Weingarten surface such that the Gauss curvature $K$ and mean curvatures $H$ satisfy the linear Weingarten condition
\begin{equation}
    aK+2bH+c=0,
\end{equation}
then its parallel surface $M^{\lambda}$ is a linear Weingarten surface such that the Gauss curvature $\bar{K}$ and mean curvature $\bar{H}$ hold
\begin{equation}
    (a+2\lambda b+\lambda^2 c)\bar{K} +2(b+\lambda c)\bar{H}+c=0.
\end{equation}
\end{lemma}

\section{$K^{\alpha}-$translators of parallel surfaces}
\label{Sec:2}

In this section, we investigate $K^{\alpha}-$translators of parallel surfaces.

\begin{theorem}
\label{transparallel}
Let $M$ be a $K^{\alpha}-$translator along $w \in \mathbb{R}^3$.
Then, its parallel surface $M^{\lambda}$
is a $\bar{K}^{\alpha}-$translator along $w \in \mathbb{R}^3$ if and only if
we have the following relations:
\begin{itemize}
\item [i.] $\lambda K-2H=0$ or $\lambda\bar{K}+2\bar{H}=0$ for $\alpha\in\mathbb{R}-\{0\}$,

\item [ii.]  $\lambda^2 K-2\lambda H+2=0$ or $\lambda^2\bar{K}+2\lambda\bar{H}+2=0$ for odd numbers $\alpha$,
\end{itemize}
where $\lambda\in\mathbb{R}$.
Thus, the surface $M$ and its parallel surface $M^{\lambda}$ move under the same $K^{\alpha}$--flow along the same direction $w$.
\end{theorem}

\begin{proof}
Let $M$ be a $K^{\alpha}-$translator along $w \in \mathbb{R}^3$. Then, it satisfies the equation \eqref{deftranslator}.
Since $M^\lambda$ is a parallel surface to $M$, we have $\bar{U}=\varepsilon U$ and \eqref{Gaussparallel}.
Thus, the equation \eqref{deftranslator} can be rewrite
\begin{equation}
\langle \bar{U}, w \rangle=\varepsilon\bar{K}^{\alpha}(1-2\lambda H+ \lambda^{2} K)^{\alpha}.
\end{equation}
If $M^{\lambda}$ is a $\bar{K}^{\alpha}-$translator
along $w \in \mathbb{R}^3$, then satisfies $(1-2\lambda H+ \lambda^{2} K)^\alpha=\varepsilon$. Thus, we have the following two cases:

\textit{Case (i).} $1-2\lambda H+ \lambda^{2} K >0$,
that is, $\varepsilon=1$. For $\alpha\in\mathbb{R}-\{0\}$,
we have $1-2\lambda H+ \lambda^{2} K=1$ which gives $\lambda K- 2 H=0$.
Considering this equation, \eqref{Gaussparallel} and \eqref{meancurvparallel} yield $\bar{K}=K$ and $\bar{H}=-H$.
Thus, $\lambda\bar{K}+2\bar{H}=0$ is also satisfied. Also, it can be seen that the parallel surface $M^{\lambda}$ is a $K^{\alpha}-$translator
along $w \in \mathbb{R}^3$.

\textit{Case (ii).} $1-2\lambda H+ \lambda^{2} K <0$, that is, $\varepsilon=-1$.
For $\alpha \in \mathbb{Z}-2\mathbb{Z}$, $1-2\lambda H+ \lambda^{2} K=-1$ holds. Thus, we have $\lambda^2 K-2\lambda H+2=0$.
Considering this equation, \eqref{Gaussparallel} and \eqref{meancurvparallel} yield $\bar{K}=-K$ and $\bar{H}=-H-\lambda\bar{K}$.
Thus, $\lambda^2\bar{K}+2\lambda\bar{H}+2=0$ is also satisfied. Similarly, the parallel surface $M^{\lambda}$ is a $K^{\alpha}-$translator
along $w \in \mathbb{R}^3$.

Conversely, let $\alpha\in\mathbb{R}-\{0\}$ and $\lambda K- 2 H=0$. Then, we have  $1-2\lambda H+ \lambda^{2} K =1$ and $\varepsilon=1$. In that case, the Gauss curvature and unit normal vector of the parallel surface $M^{\lambda}$ are obtained as $\bar{K}=K$ and $\bar{U}=U$, respectively.
Since $M$ is a $K^{\alpha}-$translator along $w \in \mathbb{R}^3$ for $\alpha\in\mathbb{R}-\{0\}$,
its parallel surface $M^{\lambda}$ is also a $K^{\alpha}-$translator along $w \in \mathbb{R}^3$.
Let $\alpha$ be an odd number and $\lambda^2 K-2\lambda H+2=0$. Then, we get $1-2\lambda H+ \lambda^{2} K =-1$ and $\varepsilon=-1$.  In that case, the Gauss curvature and unit normal vector of the parallel surface $M^{\lambda}$ are obtained as $\bar{K}=-K$ and $\bar{U}=-U$, respectively. Since $M$ is a $K^{\alpha}-$translator along $w \in \mathbb{R}^3$, we know that $\langle U, w \rangle=K^{\alpha}$. The last equation implies that
\begin{equation}
\langle -\bar{U}, w \rangle=\bar{K}^{\alpha}(-1)^{\alpha}.
\end{equation}
Thus, the parallel surface $M^{\lambda}$ becomes a $K^{\alpha}-$translator along $w \in \mathbb{R}^3$.
\end{proof}

From Theorem \ref{transparallel}, we can give the following corollary without proof.

\begin{corollary}
\label{corparallel1}
Let $M$ be a $K^{\alpha}-$translator along $w \in \mathbb{R}^3$.
If its parallel surface $M^{\lambda}$ is a $\bar{K}^{\alpha}-$translator along $w \in \mathbb{R}^3$,
then $M$ is a $(K,H)$--linear Weingarten surface.
Also, its parallel surface $M^{\lambda}$ is a $(K,H)$--linear Weingarten surface.
\end{corollary}

\begin{theorem}
Let $M$ be a $K^{\alpha}-$translator along $w \in \mathbb{R}^3$ and $1-2\lambda H+ \lambda^{2} K >0$ for $\lambda \in \mathbb{R}$.
Then, its parallel surface $M^{\lambda}$ is a $\bar{K}^{\alpha}-$translator along $w \in \mathbb{R}^3$
if and only if the surface $M^{\frac{\lambda}{2}}$ is minimal.
\end{theorem}

\begin{proof}
We assume that $M$ is a $K^{\alpha}-$translator along $w \in \mathbb{R}^3$ such that its Gauss and mean curvatures satisfy $1-2\lambda H+ \lambda^{2} K >0$ for $\lambda \in \mathbb{R}$. From Theorem \ref{transparallel}, if
$M^{\lambda}$ is a $\bar{K}^{\alpha}-$translator then we get $2H-\lambda K=0$.
On the other hand, by replacing $\lambda$ with $\frac{\lambda}{2}$ in the equation \eqref{meancurvparallel}, the mean curvature $\bar{\bar{H}}$ of $M^{\frac{\lambda}{2}}$ is given by
\begin{equation}
\bar{\bar{H}}=\frac{H-\frac{\lambda}{2} K}{1-2(\frac{\lambda}{2})H+ (\frac{\lambda}{2})^{2}K},
\end{equation}
and so $\bar{\bar{H}}=0$.

Conversely, if $\bar{\bar{H}}=0$, then $2H-\lambda K=0$. Thus, $M^{\lambda}$ is a $\bar{K}^{\alpha}-$translator along $w \in \mathbb{R}^3$.

Thus, the proof is completed.
\end{proof}

\begin{theorem}
Let $M$ be a $K^{\alpha}-$translator along $w \in \mathbb{R}^3$
and $1-2\lambda H+ \lambda^{2} K <0$ for odd numbers $\alpha$ and $\lambda \in \mathbb{R}$.
Then, its parallel surface $M^{\lambda}$ of $M$ is a $\bar{K}^{\alpha}-$translator along $w \in \mathbb{R}^3$
if and only if the surface $M^{\frac{\lambda}{2}}$ is a surface with constant Gauss curvature $\bar{\bar{K}}=-\frac{4}{\lambda^2}$ .
\end{theorem}

\begin{proof}
We assume that $M$ is a $K^{\alpha}-$translator along $w \in \mathbb{R}^3$ such that
its Gauss and mean curvatures satisfy $1-2\lambda H+ \lambda^{2} K <0$ for odd numbers $\alpha$.
From Theorem \ref{transparallel},
if $M^{\lambda}$ is a $\bar{K}^{\alpha}-$translator, then we get $\lambda^2 K-2\lambda H+2=0$.
On the other hand, by replacing $\lambda$ with $\frac{\lambda}{2}$ in the equation \eqref{Gaussparallel},
the Gauss curvature $\bar{\bar{K}}$ of $M^{\frac{\lambda}{2}}$ is given by
\begin{equation}
\label{Gausslambda2}
\bar{\bar{K}}=\frac{K}{1-\lambda H+\frac{\lambda^2}{4}K}.
\end{equation}
Substituting $\lambda^2K=-2(1-\lambda H)$ in \eqref{Gausslambda2}, we obtain the desired result.

The converse of the proof can be obtained directly.
\end{proof}

\begin{theorem}
\label{generalWeingarten}
Let $M$ be a $(K,H)$--linear Weingarten surface with $aK+2bH+c=0$ and $b^2-ac\neq 0$ for $a, b\neq 0, c \in \mathbb{R}$.
If $M$ is a $K^{\alpha}-$translator along $w \in \mathbb{R}^3$,
then its parallel surface $M^\lambda$ is a $\bar{K}^{\alpha}-$translator along $\bar{w}={\epsilon}{\mu^{\alpha}}w$,
where $\mu=\frac{b^2}{b^2-ac}$ and $\lambda=-\frac{a}{b}$ .
\end{theorem}

\begin{proof}
Assume that $M$ is a linear Weingarten surface with $aK+2bH+c=0$ and $b^2-ac\neq 0$ for $a, b\neq 0, c \in \mathbb{R}$.
The Gauss curvature of the parallel surface $M^\lambda$ for $\lambda=-\frac{a}{b}$ is computed by
\begin{equation}
\bar{K}  =  \frac{K}{1+\frac{2a}{b}H+\frac{a^2}{b^2}K}.
\end{equation}
Using $aK+2bH+c=0$, we obtain $\bar{K} =\frac{b^2}{b^2-ac}K$. Say $\mu=\frac{b^2}{b^2-ac}$.
Since $M$ is a $K^{\alpha}-$translator along $w \in \mathbb{R}^3$, it holds $\langle U, w \rangle=K^{\alpha}$.
Then, we have
\begin{align}
\langle \varepsilon \bar{U}, w \rangle & =  \frac{1}{\mu^\alpha}\bar{K}^\alpha, \\ \nonumber
\langle \bar{U}, \varepsilon{\mu^{\alpha}}w\rangle & =  \bar{K}^\alpha.
\end{align}
Thus, the parallel surface $M^{-\frac{a}{b}}$ is a $\bar{K}^{\alpha}-$translator along $\bar{w}={\epsilon}{\mu^{\alpha}}w$. The proof is completed.
\end{proof}

From Theorem \ref{generalWeingarten}, we can say that if the surface $M$ is $K^{\alpha}$--translator along $w\in\mathbb{R}^3$
and it is a $(K,H)$--linear Weingarten surface, then its parallel surface $M^{\lambda}$ is also $K^{\alpha}$--translator
along $\bar{w}={\epsilon}{\mu^{\alpha}}w\in\mathbb{R}$.

\begin{remark}
Theorem \ref{transparallel} implies that there does not exist $K^{\alpha}$--translators on the parallel surface of surface with constant mean curvature,
including minimal surfaces.
\end{remark}

\section{$K^{\alpha}-$translators of canal surfaces}
\label{Sec:3}
In this section, we first give the parametrizations of canal surfaces. Then, we study $K^{\alpha}$--translators of canal surfaces.

A canal surface $M$ can be parametrized as follows
\begin{equation}
\label{positioncanal1}
X(s,\theta)=c(s)+r(s)(\sqrt{1-r'^2(s)}\cos{\theta}N+\sqrt{1-r'^2(s)}\sin{\theta}B-r'(s)T),
\end{equation}
where $c(s)$ is a unit speed curve parametrized by arc-length $s$,
$\{T,N,B\}$ is a Frenet frame of $c(s)$.
In the sequel, $T,N,B$ is the unit tangent, principal normal and binormal vector fields, respectively.
The curve $c(s)$ is called the spine or center curve
and $r(s)$ is a radial function of $M$.
From \eqref{positioncanal1}, we may assume $r'(s)=-\cos{\varphi}$ for some smooth function $\varphi=\varphi(s)$.
Then, the canal surface $M$ can be written as
\begin{equation}
\label{positioncanal2}
X(s,\theta)=c(s)+r(s)(\sin{\varphi}\cos{\theta}N+\sin{\varphi}\sin{\theta}B+\cos{\varphi}T),
\end{equation}
where $s\in [0,\ell], \theta\in [0, 2\pi), \varphi\in [0,\pi)$ and $\ell$ is the total length of $c(s)$.
The Frenet equations of a regular space curve are given by
\begin{equation}
    \label{Freneteq}
    T'=\kappa N, \;\; N'=-\kappa T+\tau B,\;\; B'=-\tau N,
\end{equation}
where the prime denotes the differentiation with respect to $s$.
The function $\kappa$ and $\tau$ are called the curvature and torsion, respectively.

In particular, if the spine curve $c(s)$ is a straight line, then the canal surfaces become surfaces of revolution.
If the radial function $r(s)$ is a constant, then the canal surfaces are called by tubes.
Thus, all the surfaces of revolution and tubes are subclass of the canal surfaces.

The unit normal vector field to $M$ is given by
\begin{equation}
\label{normal}
    U=\cos{\varphi} T + \sin{\varphi}\cos{\theta} N + \sin{\varphi}\sin{\theta} B.
\end{equation}
Notice that the equation \eqref{positioncanal2} is written by $X(s,\theta)=c(s)+r(s) U$.
The Gaussian curvature $K$ and the mean curvature $H$ of $M$ are given by, respectively
\begin{align}
    \label{Gausscurv}
    K&=\frac{Q}{rP},\\
    \label{meancurv}
    H&=-\frac{2P+\sin^2{\varphi}}{2rP},
\end{align}
where $P=rr''+r\kappa\sin{\varphi}\cos{\theta}-\sin^2{\varphi}$ and $Q=r''+\kappa\sin{\varphi}\cos{\theta}$.
\begin{theorem}
\cite{YHKim}
The mean curvature $H$ and the Gauss curvature $K$ of a canal surface given by \eqref{positioncanal2} hold
    \begin{equation}
    \label{linear}
        H=-\frac{1}{2}(Kr+\frac{1}{r}).
    \end{equation}
\end{theorem}

\begin{proposition}
\label{propcanal}
Let $M$ be a $K^{\alpha}$--translator with speed $w$. If $M$ is a canal surface, then the tangent vector of the spine
or center curve of $M$ is parallel to $w$.
\end{proposition}

\begin{proof}
Assume that a canal surface $M$ is a $K^{\alpha}$--translator with speed $w$ and $w$ can be expressed as $w=w_1 T + w_2 N + w_3 B$.
Considering \eqref{deftranslator} and \eqref{normal}, we obtain
\begin{equation}
K^{\alpha}= w_1 \cos{\varphi} + w_2 \sin{\varphi} \cos{\theta} + w_3 \sin{\varphi} \sin{\theta}.
\end{equation}
Since $\{1, \cos{\theta}, \sin{\theta}\}$ is linearly independent set,
we get the following equations:
\begin{align}
\label{transeq1}
    &-K^{\alpha}+\cos{\varphi} w_1=0,\\
\label{transeq2}
    &\sin{\varphi} w_2=\sin{\varphi} w_3=0.
\end{align}
The equation \eqref{transeq2} gives two cases as $\sin{\varphi}=0$ or $w_2=w_3=0$.
If $\sin{\varphi}=0$, then we have $r'(s)=-1$, that is, $Q=0$.
From \eqref{Gausscurv}, it can be easily seen that $K=0$ which can not be.
Thus, $w=w_1 T$ which means that the speed $w$ is parallel to $T$.
\end{proof}

\begin{theorem}
\label{thmcanaltrans}
Let $M$ be a $K^{\alpha}$--translator with the speed $w$.
If $M$ is a canal surface, then the spine curve or the center curve of $M$ is a straight line.
Thus, the canal surface $M$ becomes a surface of revolution parametrized by
\begin{equation}
\label{positionrot}
    X(s,\theta)=(r(s)\sin{\varphi(s)}\cos{\theta}, r(s)\sin{\varphi(s)}\sin{\theta}, r(s)\cos{\varphi(s)}+s),
\end{equation}
where the radius function $r(s)$ satisfies
\begin{equation}
\label{eqforr}
    \left(\frac{r''}{r(rr''-1+r'^2)}\right)^{\alpha}=-r'.
\end{equation}
and $w$ equals the tangent vector of the spine curve or the center curve.
\end{theorem}

\begin{proof}
Assume that a canal surface $M$ is a $K^{\alpha}$--translator with the speed $w$. Then, from Proposition \ref{propcanal},
we know that the vector $w$ is parallel to the tangent vector of
the spine or center curve $c(s)$.
From the equations \eqref{Gausscurv} and \eqref{transeq1}, we have
\begin{equation}
\label{eq1}
(1-r^2(w_1 \cos{\varphi})^{\frac{1}{\alpha}})r''+ (w_1 \cos{\varphi})^{\frac{1}{\alpha}}\sin^2{\varphi}\; r
+ \kappa\sin{\varphi}( 1-r^2 (w_1\cos{\varphi})^{\frac{1}{\alpha}})\cos{\theta}=0.
\end{equation}
Since $r$ is a function of $s$, $\kappa\sin{\varphi}( 1-r^2(w_1 \cos{\varphi})^{\frac{1}{\alpha}})=0$.
We know that $\sin{\varphi}\neq 0$ due to the fact that $K\neq 0$.
On the other hand,
if $1-r^2(w_1\cos{\varphi})^{\frac{1}{\alpha}}=0$,
then the equation \eqref{eq1} becomes $(w_1 \cos{\varphi})^{\frac{1}{\alpha}}\sin^2{\varphi}\; r=0$
which implies $\cos{\varphi}=0$. In this case, $r'=-\cos{\varphi}=0$. Thus,
a canal surface $M$ becomes a tube surface whose Gaussian curvature $K$ is zero.
Hence, $\kappa$ must be zero which says that the curve $c(s)$ is a straight line.
Moreover, the canal surface $M$ becomes a surface of revolution.

Now, we will determine parametrization such a surface of revolution.
Without loss of generality, we may assume the spine curve is $c(s)=(0,0,s)$.
Then, we have $T=(0,0,1), \; N=(1,0,0),\; B=(0,1,0)$.
From the equation \eqref{positioncanal2}, we obtain the position vector
of a surface of revolution given by \eqref{positionrot}.
Also, we can choose $w_1=1$, that is, $w=T=(0,0,1)$.
In this case, the equation \eqref{Gausscurv} becomes $\displaystyle{K=\frac{r''}{r(rr''-1+r'^2)}}$.
Combining these with \eqref{transeq1}, we obtain the equation in \eqref{eqforr}.
\end{proof}

From the proof of Theorem \ref{thmcanaltrans}, we give the following corollary.
\begin{corollary}
\label{cortubetrans}
There do not exist $K^{\alpha}$--translator with any speed $w$ on a tube surface.
\end{corollary}

Considering Theorem \ref{thmcanaltrans} and Corollary \ref{cortubetrans}, we have the following remark.

\begin{remark}
\label{corcanaltrans}
A canal surface, except a surface of revolution, does not become $K^{\alpha}$--translator with any speed $w$.
\end{remark}

From now on, we consider $K^{\alpha}$--translators on a surface of revolution given by \eqref{positionrot}.

Although we are not able to find the general solution of the differential equation in \eqref{eqforr} for any constant $\alpha$,
by a direct calculation, we will give the solutions and examples of the differential equation \eqref{eqforr} for $\alpha=1$ and $\alpha=-1/2$ as follows.

\begin{corollary}
\label{alpha1}
Let $M$ be a surface of revolution parametrized by \eqref{positionrot}. Then,
$M$ is a $K$--translator with the speed $w=(0,0,1)$ where the radius function $r(s)$ is given as the solution of the following integral
\begin{equation}
\label{alpha11}
s=c_2+\int\frac{r^2}{-1\pm \sqrt{r^4+c_1 r^2+1}}dr,
\end{equation}
where $c_1$ and $c_2$ are constants.
\end{corollary}

\begin{corollary}
\label{alpha12}
Let $M$ be a surface of revolution parametrized by \eqref{positionrot}. Then,
$M$ is a $K^{-1/2}$--translator with the speed $w=(0,0,1)$ where the radius function $r(s)$ is given as the solution of the following integrals:
\begin{equation}
\label{alpha121}
s=c_2\pm\int\frac{dr}{\sqrt{r^2-\sqrt{c_1+r^4-2r^2}}}
\end{equation}
or
\begin{equation}
\label{alpha122}
s=c_2\pm\int\frac{dr}{\sqrt{\sqrt{r^4-2r^2+c_1}+r^2}},
\end{equation}
where $c_1$ and $c_2$ are constants.
\end{corollary}

\begin{example}
For $c_1=-2$, the integral given in \eqref{alpha11} occurs
\begin{equation}
\label{radius1}
r(s)-\sqrt{2}\tanh^{-1}\left(\frac{r(s)}{\sqrt{2}}\right)+c_2=s.
\end{equation}
Then, the surface of revolution $M$ given by \eqref{positionrot} with the radius function $r(s)$ defined in \eqref{radius1} explicitly
is given by
\begin{equation}
x(s,\theta)=\left(\frac{2\sqrt{r^2(s)-1}}{r(s)}\cos{\theta}, \frac{2\sqrt{r^2(s)-1}}{r(s)}\sin{\theta},
\frac{2-r^2(s)}{r(s)}+s\right),
\end{equation}
where $s\in [0, 0.7539]$ and $\theta\in [0, 2\pi)$.

For $c_1=1$ and $r^2-1<0$, the integral given in \eqref{alpha121} becomes
\begin{equation}
\label{radius2}
\pm \frac{1}{\sqrt{2}}\ln \left({r(s)+\sqrt{r^2(s)-\frac{1}{2}}}\right)+c_2=s.
\end{equation}
Then, the surface of revolution $M$ given by \eqref{positionrot} with the radius function $r(s)$ defined in \eqref{radius2} explicitly
is given by
\begin{equation}
x(s,\theta)=\left(r(s)\sqrt{2-2r^2(s)}\cos{\theta}, r(s)\sqrt{2-2r^2(s)}\sin{\theta},
-r(s)\sqrt{2-2r^2(s)}+s\right),
\end{equation}
where $s\in [0.20, 0.62]$ and $\theta\in [0, 2\pi)$.
\begin{figure}[htbp]
\includegraphics[height=40mm]{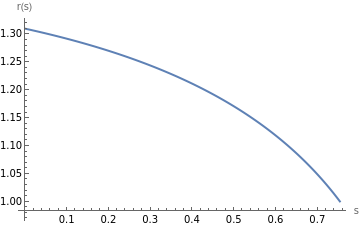}
\label{fig1}
\includegraphics[height=40mm]{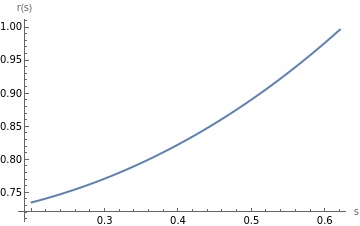}
\caption {The radius function $r(s)$ of rotational $K^{\alpha}$ surfaces: $\alpha=1$(left) and $\alpha=-1/2$ (right)}
\label{fig3}
\end{figure}
\begin{figure}[htbp]
\includegraphics[height=55mm]{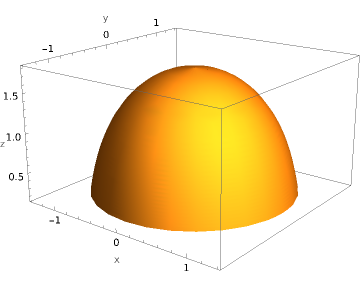}
\label{fig2}
\includegraphics[height=60mm]{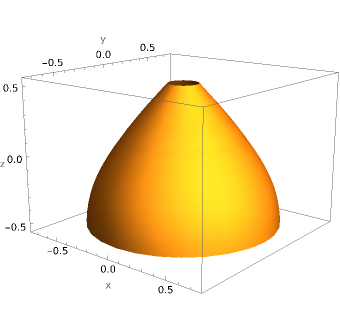}
\label{fig4}
\caption {Rotational $K^{\alpha}$ surfaces: $\alpha=1$(left) and $\alpha=-1/2$ (right)}
\label{fig2}
\end{figure}
\end{example}

Now, we consider $K^{\alpha}$--translators on $(K,H)$--linear Weingarten canal surfaces in $\mathbb{E}^3$.
\begin{theorem}
\label{linearWeingartencanal1}
\cite{YHKim}
A canal surface $M$ is a $(K,H)$--linear Weingarten canal surface such that $aK+bH=1$ if and only if it is one of the following surfaces:
\begin{itemize}
\item a tube with radius $r=-\frac{b}{2}$,

\item a surface of revolution such as
\begin{equation*}
    x(s,\theta)=(r(s)\sin{\varphi(s)}\cos{\theta}, r(s)\sin{\varphi(s)}\sin{\theta}, r(s)\cos{\varphi(s)}+s),
\end{equation*}
where $r(s)$ is given by
\begin{equation}
\label{solvingWC2}
s=c_2\pm\int \sqrt{\frac{r^2+br-a}{r^2+br-a-c_1}} dr.
\end{equation}
\end{itemize}
\end{theorem}

\begin{theorem}
\label{linearWeingartencanal2}
A canal surface $M$ given by \eqref{positioncanal2} is a $(K,H)$--linear Weingarten canal surface such that $H=cK$
if and only if it is a surface of revolution given by
\begin{equation*}
    X(s,\theta)=(r(s)\sin{\varphi(s)}\cos{\theta}, r(s)\sin{\varphi(s)}\sin{\theta}, r(s)\cos{\varphi(s)}+s),
\end{equation*}
where $r(s)$ satisfies
\begin{equation}
\label{solvingWC1}
s=c_2\pm\int \sqrt{\frac{r+c}{r+c+c_1}} dr
\end{equation}
for arbitrary constants $c, c_1$ and $c_2$.
\end{theorem}

\begin{proof}
We suppose that the mean curvature and Gaussian curvature of the canal surface $M$ given by \eqref{positioncanal2} satisfy $H=cK$. By using \eqref{linear} we get
\begin{equation}
  \label{Koflw2}
      K=-\frac{1}{r(r+2c)}.
\end{equation}
If we consider \eqref{Gausscurv} and \eqref{Koflw2} together, we have
\begin{equation}
\label{difdenk}
      2(r+c)\kappa\sin{\varphi}\cos{\theta}+2r''(r+c)-\sin^2{\varphi}=0.
\end{equation}
Then, from \eqref{difdenk} we obtain
\begin{equation}
\label{difdenk1}
    2(r+c)\kappa\sin{\varphi}=0
\end{equation}
and
\begin{equation}
\label{difdenk2}
    2r''(r+c)-(1-r'^2)=0.
\end{equation}
Since M is a regular surface, $\sin{\varphi}\neq0$.
Thus, from \eqref{difdenk1}, we get $r+c=0$ or $\kappa=0$.
For $r+c=0$, \eqref{difdenk2} is not provided. Thus,
$\kappa=0$ and M becomes a revolution surface.
Solving the differential equation given by \eqref{difdenk2}, we get the equality given by \eqref{solvingWC1}.
\end{proof}

\begin{theorem}
\label{Weingartentrans}
There do not exist ${K^\alpha}$--translators on a $(K,H)$--linear Weingarten surface of revolution given by \eqref{positionrot}.
\end{theorem}
\begin{proof}
    Assume that the surface $M$ given by \eqref{positionrot} is a $(K,H)$--linear Weingarten rotational surface
    which is a ${K^\alpha}$-- translator. Then, from \eqref{eqforr}, the following equation satisfies
    \begin{equation}
    \label{Ktrasnlator}
       {K^\alpha}=-r',
    \end{equation}
    where $K$ is the Gauss curvature of $M$ and depends only on $s$.
    If we take the derivative of \eqref{Ktrasnlator} with respect to $s$ and use \eqref{Ktrasnlator}, we get the following equation
    \begin{eqnarray}
    \label{alphaKr}
       \alpha\frac{K'}{K} = \frac{r''}{r'}.
    \end{eqnarray}
Then, we consider the following two cases:\\
\textit{Case (i.)} Let $M$ be a $(K,H)$--linear Weingarten rotational surface such that
\begin{eqnarray}
\label{lw1}
aK+bH=1
\end{eqnarray}
for some constants $a$ and $b\neq 0$.
By using \eqref{lw1} and \eqref{linear}, we have the Gauss curvature $K$ of $M$ as:
 \begin{eqnarray}
 \label{Koflw1}
     {K}=\frac{b+2r}{2ar-br^2}.
 \end{eqnarray}
By differentiating \eqref{Koflw1}, we get
\begin{equation}
  \label{Koflw1derv}
     K'=\frac{2br'(r^2+br-a)}{(2ar-br^2)^2}.
\end{equation}
 On the other hand, from the equation \eqref{solvingWC2} in Theorem \ref{linearWeingartencanal1}, we have
 \begin{equation}
 \label{r1derv1}
{r'} = \sqrt{\frac{ r^2+br-a-c_1}{r^2+br-a}}.
\end{equation}
for some constants $a, b$ and $c_1$.
Similarly the equation \eqref{r1derv1} implies
\begin{equation}
\label{r2derv1}
  {r''}  = \frac{ c_1(2r+b)}{(r^2+br-a)^2}.
\end{equation}
By substituting \eqref{Koflw1}, \eqref{Koflw1derv}, \eqref{r1derv1} and \eqref{r2derv1} into \eqref{alphaKr},
we find $6$th degree polynomial of
a function $r(s)$ as follows
\begin{equation}
\label{r2}
P_6r^6+P_5r^5+P_4r^4+P_3r^3+P_2r^2+P_1r+P_0=0.
\end{equation}
where
$P_6=4\alpha b, P_5=12\alpha b^2, P_4=4(\alpha b(-3a+3b^2- c_1)+bc_1), P_3=4\big(\alpha b^2(b^2-6a-2c_1)-c_1(2a-b^2)\big),
P_2=4\alpha b(3a^2-3ab^2+2ac_1-b^2c_1)-c_1(8ab-b^3), P_1=4\alpha b^2(3a^2+2ac_1)-2ab^2c_1$ and $P_0=-4\alpha b a^2(a+c_1)$.
If the polynomial of $r$ has a solution, then $r$ becomes a constant.
In this case, the surface $M$ is a tube surface and its Gaussian curvature is zero.
Thus, the polynomial holds for every $r$ which it implies that $b=0$.
Then, from \eqref{lw1}, $K$ must be a constant. It is a contradiction.

\textit{Case (ii.)} Let $M$ be a $(K,H)$--linear Weingarten canal surface such that $H=cK$ for a constant $c$.
Then, the Gauss curvature $K$ of $M$ is given by \eqref{Koflw2}.
Taking the derivative of \eqref{Koflw2}, we get
\begin{equation}
  \label{Koflw2derv}
     K'=\frac{2{r'} (c+r)}{(r^2+2cr)^2}.
\end{equation}
 On the other hand, from \eqref{solvingWC1} we have
 \begin{equation}
 \label{r1derv2}
{r'} = \sqrt{\frac{ r+c+c_1}{r+c}}.
\end{equation}
for some constants $c$ and $c_1$. Also, the derivative of \eqref{r1derv2} is obtained as:
\begin{equation}
\label{r2derv2}
  {r''}  = -\frac{ c_1}{2(r+c)^2}.
\end{equation}
Similar as Case (i.), substituting \eqref{Koflw2}, \eqref{Koflw1derv}, \eqref{r1derv2} and \eqref{r2derv2} into \eqref{alphaKr}, we get
the third degree polynomial of $r(s)$ as follows
 \begin{equation}
    \label{r2}
        4\alpha r^3+ (4\alpha(3c+c_1)-c_1)r^2+(4\alpha (3c^2+2cc_1)-2cc_1)r+4\alpha c^2(c+c_1)=0.
    \end{equation}
If there exists a solution of the polynomial \eqref{r2}, then $r$ is a constant.
In that case, the Gauss curvature of $M$ vanishes and it is a contradiction.
Thus, the polynomial satisfies for every $r$ which implies that $\alpha=0$. That is a contradiction.

Hence, the surface $M$ can not be a $K^{\alpha}$--translator and $(K,H)$--linear Weingarten surface at the same time.
\end{proof}

From the proof of Theorem \ref{Weingartentrans}, it can been seen that
a surface of revolution given by \eqref{positionrot} with constant mean curvature and minimal rotational surface
are obtained by using Case (i.) and Case (ii.) for $a=0$ and $c=0$,
respectively. Thus, we give the following corollary.

\begin{corollary}
There do not exist ${K^\alpha}$--translators on a surface of revolution given by
\eqref{positionrot} with constant mean curvature, including minimal surfaces.
\end{corollary}

Also, we give the following corollary about $K^{\alpha}$--translators on parallel surfaces of rotational surfaces
obtained from canal surfaces.

\begin{corollary}
     Let a surface of revolution $M$ given by \eqref{positionrot} be a $K^\alpha$--translator with the speed $w=(0,0,1)$.
     Then, there are no parallel surfaces $M^{\lambda}$ which is $\bar{K}^\alpha$--translator with same speed $w$.
\end{corollary}
\begin{proof}
    Let the surface of revolution $M$ given by \eqref{positionrot} be a $K^\alpha$--translator with the speed $w=(0,0,1)$.
    From Theorem \ref{transparallel}, if the parallel surfaces $M^{\lambda}$ of $M$ are $\bar{K}^\alpha$--translator with same direction $w$,
    then M becomes $(K,H)$--linear Weingarten such that $\lambda K-2H=0$ or  $\lambda^2 K-2\lambda H+2=0$.
    However, from Theorem \ref{Weingartentrans}, we show that there do not exist ${K^\alpha}$--translators on any $(K,H)$--linear Weingarten rotational surface. Thus, it completes the proof.
\end{proof}

\section{Conclusion}
In this paper, we investigate $K^{\alpha}$--translators on offset surfaces (parallel surfaces and canal surfaces) in 3-dimensional Euclidean space $\mathbb{E}^3$. First, we find out the condition for two parallel surfaces to become $K^{\alpha}$--translators with same speed $w$. Then, we examine $K^{\alpha}$--translators on canal surfaces and we show that if a canal surface is $K^{\alpha}$--translator, then it becomes a surface of revolution. Moreover, we give some examples for moving a surface of revolution under $K$--flow and $K^{-1/2}$--flow along a direction $w=(0,0,1)$ and we draw the graphs such surfaces by using Wolfram Mathematica 10.4. Finally, we prove that there does not exist $K^{\alpha}$--translators on the parallel surfaces of rotational surfaces formed
from canal surfaces.

In the future, we will study on $K^{\alpha}-$translators of offset surfaces in Minkowski 3--space.



\bibliographystyle{elsarticle-num}
\biboptions{compress}

\end{document}